\documentclass[british]{amsart}
\usepackage[T1]{fontenc}
\usepackage[latin9]{inputenc}
\usepackage[a4paper]{geometry}
\usepackage{amstext}
\usepackage{amsthm}
\usepackage{amssymb}
\usepackage{setspace}
\makeatletter
\numberwithin{equation}{section}
\numberwithin{figure}{section}
\theoremstyle{plain}
\newtheorem{thm}{\protect\theoremname}[section]
\theoremstyle{plain}
\newtheorem{cor}[thm]{\protect\corollaryname}
\theoremstyle{remark}
\newtheorem{rem}[thm]{\protect\remarkname}

\usepackage{amscd}
\usepackage{mathptmx}
\usepackage{bbm}
\usepackage{paralist}

\newcommand{\e}{\mathrm{e}}

\newcommand{\N}{\mathbb{N}}

\newcommand{\R}{\mathbb{R}}

\renewcommand{\Pi}{\pi}
\renewcommand{\emptyset}{\varnothing}

\newcommand\diam{\mathrm{diam}}

\DeclareMathOperator*{\Hol}{\textnormal{H\"ol}}

\usepackage{hyperref}
\@ifundefined{textmu}
 {\usepackage{textcomp}}{}

\makeatother

\usepackage{babel}
\providecommand{\corollaryname}{Corollary}
\providecommand{\remarkname}{Remark}
\providecommand{\theoremname}{Theorem}

\begin{document}
\title{Multifractal Formalism for Generalised local dimension spectra of
Gibbs measures on the real line}
\begin{abstract}
We extend the multifractal formalism for the local dimension spectrum
of a Gibbs measure $\mu$ supported on the attractor $\Lambda$ of
a conformal iterated functions system on the real line. Namely, for
$\alpha\in\R$, we establish the multifractal formalism for the Hausdorff
dimension of the set of $x\in\Lambda$ for which the $\mu$-measure
of a ball of radius $r_{n}$ centred at $x$ obeys a power law $r_{n}{}^{\alpha}$,
for a sequence $r_{n}\rightarrow0$. This allows us to investigate
the H\"older regularity of various fractal functions, such as distribution
functions and conjugacy maps associated with conformal iterated function
systems.
\end{abstract}

\author{Johannes Jaerisch}
\address{Department of Mathematics, Faculty of Science and Engineering, Shimane
University, Nishikawatsu 1060, Matsue, Shimane, 690-8504 Japan}
\email{jaerisch@riko.shimane-u.ac.jp}
\urladdr{http://www.math.shimane-u.ac.jp/\textasciitilde jaerisch}
\author{Hiroki Sumi}
\address{Course of Mathematical Science, Department of Human Coexistence, Graduate
School of Human and Environmental Studies, Kyoto University Yoshida-nihonmatsu-cho,
Sakyo-ku, Kyoto 606-8501, Japan }
\email{sumi@math.h.kyoto-u.ac.jp}
\urladdr{http://www.math.h.kyoto-u.ac.jp/\textasciitilde sumi/index.html}
\thanks{\today}
\maketitle

\section{Introduction and Statement of Results}

Multifractal analysis has its origin in statistical physics (see \cite{Mandelbrot:74,FrischParisi:85,Halsey:86}).
For the mathematical theory of multifractal formalism and its relation
to thermodynamic formalism, we recommend \cite{pesindimensiontheoryMR1489237}.
In this article we extend the multifractal formalism for local dimension
spectra of Gibbs measures supported on the attractor $\Lambda\subset\R$
of a conformal iterated function system on the real line. For a Borel
measure $\mu$ on $\R$ and with $B(x,r):=\left\{ y\in\R\mid\left|y-x\right|<r\right\} $,
the local dimension spectrum of $\mu$ is given by the level sets
\[
\mathcal{F}(\alpha):=\left\{ x\in\Lambda\mid\lim_{r\rightarrow0}\frac{\log\mu\left(B(x,r)\right)}{\log r}=\alpha\right\} ,\quad\alpha\in\R.
\]
It is well known (\cite{MR1479016,MR1435198,MR1439809,pesindimensiontheoryMR1489237})
for finitely generated contracting conformal iterated function systems
satisfying the open set condition that, if $\mu$ is the Gibbs measure
of a H\"older continuous potential, then the dimension spectrum $f$
given by 
\[
f(\alpha):=\dim_{H}\mathcal{F}(\alpha)
\]
is equal to the Legendre transformation of a certain function $t$
defined implicitly by solving a topological pressure equation. To
denote this state of affairs, we say that the multifractal formalism
holds. In our setting, the function $t$ will be given by equation
\eqref{eq:pressure equation} below. For related results on the multifractal
formalism for non-uniformly hyperbolic systems and graph-directed
constructions, we refer to \cite{MR2755925,MR2525196}.

For self-conformal measures of finitely generated contracting conformal
iterated function systems it is also known (\cite{barreiraschmelingMR1759398})
that the set of divergence points 
\[
\left\{ x\in\Lambda\mid\liminf_{r\rightarrow\infty}\frac{\log\mu\left(B(x,r)\right)}{\log r}<\limsup_{r\rightarrow\infty}\frac{\log\mu\left(B(x,r)\right)}{\log r}\right\} 
\]
has full Hausdorff dimension, unless the the spectrum $f$ is degenerate.
In \cite{MR2864378} it is shown that, for $x\in\Lambda$, the set
of accumulation points $A(x)$ of $\left(\log(r)^{-1}\log\mu\left(B(x,r)\right)\right)_{r\ge0}$
is either a singleton or a closed interval. Moreover, it is shown
that 
\[
\dim_{H}\left\{ A(x)=[a,b]\right\} =\inf_{\alpha\in\left[a,b\right]}\dim_{H}\mathcal{F}(\alpha).
\]
Recently, the following modified level sets have attracted a lot of
attention in studying the regularity of various fractal functions
(\cite{JS13b}, \cite{JS2018}, \cite{MR3746637}, \cite{MR3816651},
\cite{MR3771122}). For $\alpha\in\R$, we define 
\[
R_{*}(\alpha):=\left\{ x\in\Lambda\mid\liminf_{r\rightarrow0}\frac{\log\mu\left(B(x,r)\right)}{\log r}=\alpha\right\} ,\quad R^{*}(\alpha):=\left\{ x\in\Lambda\mid\limsup_{r\rightarrow0}\frac{\log\mu\left(B(x,r)\right)}{\log r}=\alpha\right\} .
\]
In \cite{JS13b} it is shown that the level sets $R_{*}(\alpha)$
and $R^{*}(\alpha)$ satisfy the multifractal formalism in the context
of semigroups of rational maps on the Riemann sphere satisfying the
separating condition. Allaart (\cite{MR3771122}) proved the multifractal
formalism for the level sets $R_{*}(\alpha)$ for self-similar measures
satisfying the open set condition. In \cite{JS2018} we established
the multifractal formalism for $R_{*}(\alpha)$ for self-conformal
measures supported on attractors of conformal iterated function systems
satisfying the open set condition. Allaart (\cite{MR3771122}) raised
the question whether the multifractal formalism holds for $R^{*}(\alpha)$.
Moreover, it is of interest whether the formalism for $R_{*}$ and
$R^{*}$ extends to arbitrary Gibbs measures of H\"older continuous
potentials. 

In this article we establish the multifractal formalism for the generalised
level sets 
\[
\mathcal{F}^{*}(\alpha):=\left\{ x\in\Lambda\mid\exists(r_{k})\rightarrow\infty\,\,\lim_{k\rightarrow\infty}\frac{\log\mu\left(B(x,r_{k})\right)}{\log r_{k}}=\alpha\right\} ,\quad\alpha\in\R,
\]
for Gibbs measures of H\"older continuous potentials supported on
attractors of conformal iterated function systems satisfying the open
set condition. Combining with the previously known results for $\mathcal{F}(\alpha)$,
we thus obtain that the multifractal formalism holds in particular
for the level sets $R_{*}(\alpha)$ and $R^{*}(\alpha)$.

Let us now introduce the necessary terminology to state our main theorem
precisely. For an index set $I$ with $\#I<\infty$ let $\Phi=(\phi_{i})_{i\in I}$,
$\phi_{i}:X\rightarrow X$, be a contracting $\mathcal{C}^{1+\epsilon}$
conformal iterated function system on a compact subset $X\subset\R$.
We refer to \cite{MR1387085} for further details and basic properties.
For simplicity, we will also assume that $\phi_{i}'>0$ on $X$ for
each $i\in I$. Let $\pi:I^{\N}\rightarrow\R$ denote the coding map
of $\Phi$ which is for $\omega=(\omega_{1},\omega_{2},\dots)\in I^{\N}$
given by 
\[
\bigcap_{n\ge1}\phi_{\omega_{1}}\circ\dots\circ\phi_{\omega_{n}}(X)=\left\{ \pi(\omega)\right\} .
\]
Let $\Lambda:=\pi(I^{\N})$. We say that $\Phi$ satisfies the open
set condition if there exists a non-empty open interval $U\subset\R$
such that $\phi_{i}(U)\subset U$, for all $i\in I$, and $\phi_{i}(U)\cap\phi_{j}(U)=\emptyset$
for all $i,j\in I$ with $i\neq j$. To utilize the symbolic thermodynamic
formalism (see \cite{bowenequilibriumMR0442989}) we will also need
the following definitions. We denote by $\sigma:I^{\N}\rightarrow I^{\N}$
the left-shift on $I^{\N}$ which becomes a compact metric space endowed
with the shift metric. Let $\varphi:I^{\N}\rightarrow\R$ be the geometric
potential of $\Phi$ given by 
\[
\varphi(\omega):=\log\phi_{\omega_{1}}'(\pi(\sigma(\omega))),\quad\omega=(\omega_{1},\omega_{2},\dots)\in I^{\N}.
\]
Note that $\varphi$ is H\"older continuous. Let $\psi:I^{\N}\rightarrow\R$
be H\"older continuous and let $\mu_{\psi}$ denote the unique Gibbs
measure for $\psi$ in the sense of Bowen (\cite{bowenequilibriumMR0442989}).
For $n\ge1$ we denote by $S_{n}\psi:=\sum_{k=0}^{n-1}\psi\circ\sigma^{k}$
the ergodic sum. The range of the multifractal spectrum is defined
by 
\[
\alpha_{-}:=\inf_{\omega\in I^{\N}}\liminf_{n\rightarrow\infty}\frac{S_{n}\psi(\omega)}{S_{n}\varphi(\omega)},\quad\alpha_{+}:=\sup_{\omega\in I^{\N}}\limsup_{n\rightarrow\infty}\frac{S_{n}\psi(\omega)}{S_{n}\varphi(\omega)}.
\]
Recall that $\alpha_{-}\le\dim_{H}\Lambda$ with equality if and only
if $\alpha_{-}=\alpha_{+}$ (see \cite{pesindimensiontheoryMR1489237}).
Since $\varphi<0$ we have that for each $\beta\in\R$ there exists
a unique $t(\beta)$ such that 
\begin{equation}
\mathcal{P}(t(\beta)\varphi+\beta\psi)=0,\label{eq:pressure equation}
\end{equation}
where $\mathcal{P}(f)$ refers to the topological pressure of a continuous
function $f:I^{\N}\rightarrow\R$ with respect to $\sigma$. We denote
by 
\[
t^{*}(\alpha):=\sup_{x\in\R}\left(\alpha x-t(x)\right)
\]
the Legendre transform of $t$. 
\begin{thm}
\label{thm:mf}Let $\Phi=(\phi_{i}:X\rightarrow X)_{i\in I}$ be a
finitely generated $\mathcal{C}^{1+\epsilon}$ conformal iterated
function system on $X\subset\R$ satisfying the open set condition.
Let $\psi:I^{\N}\rightarrow\R$ be H\"older continuous and let $\mu:=\mu_{\psi}\circ\pi^{-1}$.
Then we have for all $\alpha\in[\alpha_{-},\alpha_{+}],$
\[
\dim_{H}R_{*}(\alpha)=\dim_{H}R^{*}(\alpha)=\dim_{H}\mathcal{F}(\alpha)=\dim_{H}\mathcal{F}^{*}(\alpha)=-t^{*}(-\alpha),
\]
and for $\alpha\not\notin[\alpha_{-},\alpha_{+}]$ we have $R_{*}(\alpha)=R^{*}(\alpha)=\mathcal{F}(\alpha)=\mathcal{F}^{*}(\alpha)=\emptyset$. 
\end{thm}

We proceed with two applications of our main result to the regularity
of fractal functions.

\subsection{Distribution functions of Gibbs measures}

For a function $F:\R\rightarrow\R$ we define the\emph{ }pointwise
H\"older exponent of $F$ at $x\in\R$ by
\[
\Hol(F,x):=\sup\left\{ \alpha>0\mid\limsup_{y\rightarrow x}\frac{\left|F(y)-F(x)\right|}{\left|y-x\right|^{\alpha}}<\infty\right\} .
\]
Under the assumptions of Theorem \ref{thm:mf}, we consider the distribution
function of $\mu=\mu_{\psi}\circ\pi^{-1}$ given by 
\[
F_{\mu}:\R\rightarrow\left[0,1\right],\quad F_{\mu}(x):=\mu\left((-\infty,x]\right).
\]

It is shown in \cite[Lemma 5.1]{JS13b} that 
\begin{equation}
\Hol(F_{\mu},x)=\liminf_{r\rightarrow0}\frac{\log\mu\left(B(x,r)\right)}{\log r}.\label{eq:hoelder exponent as liminf}
\end{equation}

\begin{cor}
\label{cor:pointwise hoelder}Under the assumptions of Theorem \ref{thm:mf},
the distribution function $F_{\mu}:\R\rightarrow[0,1]$ of $\mu=\mu_{\psi}\circ\pi^{-1}$
satisfies for every $\alpha\in\left[\alpha_{-},\alpha_{+}\right]$
\[
\dim_{H}\left\{ x\in\Lambda\mid\Hol(F_{\mu},x)=\alpha\right\} =-t^{*}(-\alpha),
\]
where $t:\R\rightarrow\R$ is defined implicitly by $\mathcal{P}(t(\beta)\varphi+\beta\psi)=0$. 
\end{cor}

For results on points of non-differentiability and related properties
of $F_{\mu}$ we refer to \cite{MR2034020,MR2525939,MR3177871}. 

\subsection{Conjugacy maps between expanding piecewise $\mathcal{C}^{1+\epsilon}$
interval maps}

In this section we apply the multifractal formalism to conjugacy maps
between expanding piecewise $\mathcal{C}^{1+\epsilon}$ interval maps
as considered in \cite{MR2576266}. In fact, we slightly extend the
framework  by allowing one of the repellers of the expanding interval
maps to be a proper subset of $[0,1]$, whereas in \cite{MR2576266}
both repellers are equal to $[0,1]$. 

Let us now briefly introduce the setting. Let $f$ be an expanding
piecewise $\mathcal{C}^{1+\epsilon}$ interval map with $s\ge2$ full
branches, i.e., there exist closed intervals $J_{1},\dots,J_{s}\subset[0,1]$
with non-empty, pairwise disjoint interiors such that $f_{|J_{i}}$
has a $\mathcal{C}^{1+\epsilon}$ extension to a neighbourhood of
$J_{i}$ satisfying $f'_{|J_{i}}>1$ and $f(J_{i})=[0,1]$, for $1\le i\le s$.
We will always assume that the intervals $J_{1},\dots,J_{s}$ are
given in increasing order (i.e., $\sup J_{i}\le\inf J_{j}$ if $i<j$).
The repeller of $f$ is denoted by $\Lambda$ and the restriction
$f_{|\Lambda}:\Lambda\rightarrow\Lambda$ is conjugate to the shift
$\sigma:I^{\N}\rightarrow I^{\N}$ with $I:=\left\{ 1,\dots,s\right\} $.
The conjugacy is given by the coding map $\pi_{f}:I^{\N}\rightarrow\Lambda$
of the conformal iterated function systems $\Phi_{f}$, which is given
by the contracting inverse branches $(f_{|J_{i}})^{-1}:[0,1]\rightarrow[0,1]$,
$1\le i\le s$. Similarly, let $g:\left[0,1\right]\rightarrow\left[0,1\right]$
be an expanding $\mathcal{C}^{1+\epsilon}$ interval map with $s$
full branches. We assume that the repeller of $g$ is the interval
$\left[0,1\right]$. Again, $g$ is conjugate to $(I^{\N},\sigma)$
with conjugacy map given by the coding map $\pi_{g}:I^{\N}\rightarrow\left[0,1\right]$
of the associated conformal iterated function system $\Phi_{g}$.
Thus, there is a natural conjugacy map $\Theta:\Lambda\rightarrow\left[0,1\right]$
given by $\Theta:=\pi_{g}\circ\pi_{f}^{-1}$ satisfying $\Theta\circ f_{|\Lambda}=g\circ\Theta$
(see \cite{MR2576266} for the case $\Lambda=\left[0,1]\right]$).
Note that $\#\left(\pi_{f}^{-1}(x)\right)\le2$ and that $\pi_{g}\circ\pi_{f}^{-1}$
is a well-defined non-decreasing function on the real line, which
is strictly increasing function on $\Lambda$. In particular, we have
for $x\in\Lambda$ that 
\begin{equation}
\pi_{g}^{-1}\left((-\infty,\pi_{g}\circ\pi_{f}^{-1}(x)]\right)=\pi_{f}^{-1}\left((-\infty,x]\right).\label{eq:sets}
\end{equation}
Denote by $\varphi_{g}:I^{\N}\rightarrow\R$ (resp. $\varphi_{f}:I^{\N}\rightarrow\R)$
the geometric potential of $\Phi_{g}$ (resp. $\Phi_{f})$ given by
\[
\varphi_{g}:=-\log g'\circ\pi_{g},\quad\varphi_{f}:=-\log f'\circ\pi_{f}.
\]
Let $\lambda$ denote the Lebesgue measure on $\left[0,1\right]$
and recall that $\lambda=\mu_{\varphi_{g}}\circ\pi_{g}^{-1}$ where
$\mu_{\varphi_{g}}$ is the unique Gibbs measure for $\varphi_{g}$
on $I^{\N}$. By \eqref{eq:sets} we then have for $x\in\Lambda$,
\[
\Theta(x)=\mu_{\varphi_{g}}\circ\pi_{g}^{-1}\left((-\infty,\Theta(x)]\right)=\mu_{\varphi_{g}}\circ\pi_{g}^{-1}\left((-\infty,\pi_{g}\circ\pi_{f}^{-1}(x)]\right)=\mu_{\varphi_{g}}\circ\pi_{f}^{-1}\left((-\infty,x]\right).
\]
So, the conjugacy $\Theta:\Lambda\rightarrow[0,1]$ coincides with
the distribution function of $\mu_{\varphi_{g}}\circ\pi_{f}^{-1}$
(cf. \cite{MR2576266} for the case $\Lambda=\left[0,1\right]$).
Therefore, we have 
\[
\left\{ x\in\Lambda\mid\Hol(\Theta,x)=\alpha\right\} =\left\{ x\in\Lambda\mid\Hol(F_{\mu},x)=\alpha\right\} ,
\]
where $\mu=\mu_{\varphi_{g}}\circ\pi_{f}^{-1}$. Hence, with 
\[
\alpha_{-}=\inf_{\omega\in I^{\N}}\liminf_{n\rightarrow\infty}\frac{S_{n}\varphi_{g}(\omega)}{S_{n}\varphi_{f}(\omega)},\quad\alpha_{+}:=\sup_{\omega\in I^{\N}}\limsup_{n\rightarrow\infty}\frac{S_{n}\varphi_{g}(\omega)}{S_{n}\varphi_{f}(\omega)},
\]
we obtain the following corollary as a consequence of Corollary \ref{cor:pointwise hoelder}.
\begin{cor}
Let $f$ and $g$ be two expanding piecewise $\mathcal{C}^{1+\epsilon}$
interval maps with $s\ge2$ full branches and coding maps $\pi_{f},\pi_{g}:I^{\N}\rightarrow[0,1]$.
Let $\Lambda:=\pi_{f}(I^{\N})$ and suppose that $\pi_{g}(I^{\N})=[0,1]$.
Then the conjugacy map $\Theta:\Lambda\rightarrow\left[0,1\right]$,
given by $\Theta:=\pi_{g}\circ\pi_{f}^{-1}$, satisfies for every
$\alpha\in\left[\alpha_{-},\alpha_{+}\right]$ 
\[
\dim_{H}\left\{ x\in\Lambda\mid\Hol(\Theta,x)=\alpha\right\} =-t^{*}(-\alpha),
\]
where $t:\R\rightarrow\R$ is defined implicitly by \textbf{$\mathcal{P}(t(\beta)\varphi_{f}+\beta\varphi_{g})=0$. }
\end{cor}

\section{Proof of Theorem \ref{thm:mf}}
\begin{proof}
First observe that $\mathcal{F}(\alpha)$ is a subset of each of the
level sets considered in the theorem. Therefore, the lower bound for
the Hausdorff dimension follows from the well-known multifractal formalism
for $\mathcal{F}(\alpha)$. Therefore, to complete the proof of the
theorem, it suffices to show the upper bound for the Hausdorff dimension
of $\mathcal{F}^{*}(\alpha)$. Throughout, we may assume $I=\left\{ 1,\dots,s\right\} $,
for $s\ge2$, and $\phi _{i} (x)\leq \phi _{y}(y)$ for all  $i,j\in I$ with $i<j$ and for all $x,y\in U$,
where $U$ is the open set in the open set condition.

Denote by $I^{*}:=\bigcup_{k\ge1}I^{k}$ the set of finite words in
the alphabet $I$. We first observe that there exists $C>1$ such
that for all $\gamma\in I^{*}$ and $i\in I$, 
\[
C\cdot\diam\left(\pi\left[\gamma i\right]\right)\ge\diam\left(\pi\left[\gamma\right]\right),
\]
where $\diam(A):=\sup\left\{ |x-y|\mid x,y\in A\right\} $ refers
to the diameter of a set $A\subset\R$. 

Let $x\in\mathcal{F}^{*}(\alpha)$, $x=\pi(\omega)$ for some $\omega\in I^{\N}.$
There exists a sequence $r_{k}\rightarrow0$ such that 
\[
\lim_{k\rightarrow\infty}\frac{\log\mu\left(B(x,r_{k})\right)}{\log r_{k}}=\alpha.
\]
We define for $k\ge1$
\[
n_{k}:=\min\left\{ n\ge1\mid\diam\left(\pi\left[\omega_{1},\dots,\omega_{n}\right]\right)<Cr_{k}\right\} .
\]
For $a\in I$ and $m\in\N$ we denote $a^{m}:=(a,\dots,a)\in I^{m}$.
We define two sequences of words $(\nu_{k})$, $(\nu_{k}')\in I^{*}$,
$k\ge1$, as follows. If $(\omega_{1},\dots,\omega_{n_{k}})$ takes
the form 
\[
(\omega_{1},\dots,\omega_{n_{k}})=(\tau js^{\ell_{k}}),
\]
for some $\tau\in I^{*}$, $j\le s-1$ and $\ell_{k}\ge1$ then let
\[
\ell'_{k}:=\max\left\{ l\ge1\mid\diam\left(\pi\left[\tau(j+1)1^{l}\right]\right)\right\} \ge r_{k}.
\]
We remark that, if the constant $C$ above is large enough, then $\ell_{k}'$
is well defined. This follows because our definition of $n_{k}$ implies
$\diam\left(\pi\left[\tau j\right]\right)\ge Cr_{k}$. We then define
\[
\nu_{k}:=(\tau js^{\ell_{k}-1})=(\omega_{1},\dots,\omega_{n_{k}-1}),\quad\nu_{k}':=(\tau(j+1)1^{\ell'_{k}}).
\]
If $(\omega_{1},\dots,\omega_{n_{k}})$ takes the form $(\omega_{1},\dots,\omega_{n_{k}})=(\tau(j+1)1^{\ell_{k}}),$
then we define $\nu_{k}:=(\omega_{1},\dots,\omega_{n_{k}-1})$ and
$\nu'_{k}:=(\tau js^{\ell'_{k}})$ in a similar fashion.  Finally,
if $(\omega_{1},\dots,\omega_{n_{k}})$ satisfies $\omega_{n_{k}}\notin\{1,s\}$
then we define $\nu_{k}:=\nu_{k}':=(\omega_{1},\dots,\omega_{n_{k}-1}).$
Let $n_{k}':=|\nu_{k}'|$.

It is important to note that, by the definition of $\nu_{k}$ and
$\nu_{k}'$, we have as $k\rightarrow\infty$, 
\begin{equation}
\diam\left(\pi\left[\nu_{k}\right]\right)\asymp\diam\left(\pi\left[\nu_{k}'\right]\right)\asymp r_{k},\label{eq:comparable}
\end{equation}
where, for sequences of positive numbers $(a_{k})$ and $(b_{k})$
the notation $a_{k}\asymp b_{k}$ means that $a_{k}/b_{k}$ is bounded
away from zero and infinity.  We will show that this implies the
existence of a uniform constant $D$ such that 
\begin{equation}
\ell_{k}'\cdot\varphi(\overline{1})-D\le\ell_{k}\cdot\varphi(\overline{s})\le\ell_{k}'\cdot\varphi(\overline{1})+D,\label{eq:comparability lk}
\end{equation}
where we have set $\overline{\gamma}:=(\gamma\gamma\dots)\in I^{\N}$,
for \textbf{$\gamma\in I^{*}$.} To prove \eqref{eq:comparability lk}
suppose that $(\omega_{1},\dots,\omega_{n_{k}})=(\tau js^{\ell_{k}})$.
The other case $(\omega_{1},\dots,\omega_{n_{k}})=(\tau(j+1)1^{\ell_{k}}),$
can be handled analogously. By the bounded distortion property of
the geometric potential we have, as $k\rightarrow\infty$, 
\[
\diam\left(\pi\left[\nu_{k}\right]\right)=\diam\left(\pi\left[(\tau js^{\ell_{k}-1})\right]\right)\asymp\diam\left(\pi\left[\tau\right]\right)\e^{S_{\ell_{k}}\varphi(\overline{s})}
\]
\[
\diam\left(\pi\left[\nu_{k}'\right]\right)=\diam\left(\pi\left[\tau(j+1)1^{\ell'_{k}}\right]\right)\asymp\diam\left(\pi\left[\tau\right]\right)\e^{S_{\ell'_{k}}\varphi(\overline{1})},
\]
which proves \eqref{eq:comparability lk}.

We will show that for $x\in\Lambda=\pi(I^{\N})$ and $k\ge1$ we have
\[
B(x,r_{k})\cap\Lambda\subset\pi([\nu_{k}])\cup\pi([\nu'_{k}]).
\]

First suppose that $(\omega_{1},\dots,\omega_{n_{k}})=(\tau js^{\ell_{k}})$.
Then $\pi\left[\nu_{k}\right]=\pi\left[(\tau js^{\ell_{k}-1})\right]\supset\pi\left[(\tau js^{\ell_{k}-1}1)\right]$.
Since $x\in\pi\left[(\tau js^{\ell_{k}})\right]$ we have $x\ge\max\pi\left[(\tau js^{\ell_{k}-1}1)\right]$.
By the definition of $C$ and $n_{k}$ we have 
\[
\diam\left(\pi\left[(\tau js^{\ell_{k}-1}1)\right]\right)\ge C^{-1}\diam\left(\pi\left[(\tau js^{\ell_{k}-1})\right]\right)=C^{-1}\diam\left(\pi\left[(\omega_{1},\dots,\omega_{n_{k}-1})\right]\right)\ge r_{k}.
\]
Hence, 
\[
\pi\left[\nu_{k}\right]\supset\left[x,x-r_{k}\right]\cap\Lambda.
\]
Further, by the definition of $\ell'_{k}$ we have 
\[
\diam\left(\pi\left[\nu_{k}'\right]\right)=\diam\left(\pi\left[\tau(j+1)1^{\ell'_{k}}\right]\right)\ge r_{k},
\]
so $\left[x,x+r_{k}\right]\cap\Lambda\subset\pi\left[\nu_{k}\right]\cup\pi\left[\nu_{k}'\right]$.
This proves that $B(x,r_{k})\cap\Lambda\subset\pi([\nu_{k}])\cup\pi([\nu'_{k}])$. 

Let $\epsilon>0$. We will derive from our assumption $x\in\mathcal{F}^{*}(\alpha)$
that there exists $N\ge1$ such that for all $k\ge N$, 
\begin{equation}
\frac{S_{\left|\nu_{k}\right|}\psi(\overline{\nu_{k}})}{S_{\left|\nu_{k}\right|}\varphi(\overline{\nu_{k}})}\le\alpha+\epsilon\quad\text{or}\quad\frac{S_{\left|\nu_{k}'\right|}\psi(\overline{\nu_{k}'})}{S_{\left|\nu_{k}'\right|}\varphi(\overline{\nu_{k}'})}\le\alpha+\epsilon.\label{eq:nu or nuprime is good}
\end{equation}
To prove \eqref{eq:nu or nuprime is good}, we first note that by
the Gibbs property of $\mu_{\psi}$ we have for every $\gamma\in I^{*}$,
\[
\mu_{\psi}\left(\left[\gamma\right]\right)\le C_{\psi}\exp(S_{|\gamma|}\psi(\overline{\gamma})),
\]
where $C_{\psi}\ge0$ is a uniform constant depending on $\psi$.
Suppose for a contradiction that \eqref{eq:nu or nuprime is good}
does not hold. Then, by passing to a subsequence of $(n_{k})$ we
may assume that for all $k$ and for all $\nu\in\{\nu_{k},\nu'_{k}\}$
we have $S_{|\nu|}\psi(\overline{\nu})\big/S_{|\nu|}\varphi(\overline{\nu})>\alpha+\epsilon$,
and so, taking $C'$ large enough, we have 
\[
\mu_{\psi}\left(\left[\nu\right]\right)\le C_{\psi}\exp(S_{|\nu|}\psi(\overline{\nu}))\le C'r_{k}^{\alpha+\epsilon}.
\]
Since $B(x,r_{k})\subset\pi([\nu_{k}])\cup\pi([\nu'_{k}])$, we conclude
that 
\[
\lim_{k\rightarrow\infty}\frac{\log\mu\left(B(x,r_{k})\right)}{\log r_{k}}\ge\alpha+\epsilon.
\]
This contradiction proves \eqref{eq:nu or nuprime is good}. 

We will prove the proposition only in the case when $\alpha\le\alpha_{0}:=\int\psi\,\,d\mu_{t(0)\varphi}\big/\int\varphi\,\,d\mu_{t(0)\varphi}$.
The case $\alpha\ge\alpha_{0}$ can be considered in a similar fashion
(see also Remark \ref{rem:upper half} below). Let $\beta>0$, $\eta>0$
and let $b=t(\beta)+\beta(\alpha+\epsilon)+\eta$. We define
\[
\mathcal{C}_{\alpha+\epsilon}:=\left\{ \tau\in I^{*}\mid\frac{S_{\left|\tau\right|}\psi(\overline{\tau})}{S_{\left|\tau\right|}\varphi(\overline{\tau})}\le\alpha+\epsilon\right\} .
\]
We obtain a covering $\mathcal{C}$ of $\mathcal{F}^{*}(\alpha)$
by cylinders of sufficiently small diameters as follows. For each
$x\in\mathcal{F}^{*}(\alpha)$ we define the sequence $n_{k}$ as
above and define 
\[
\mathcal{C}:=\left\{ \nu_{k}\mid k\in\N\right\} ,
\]
where $\nu_{k}\in I^{*}$ is defined as above satisfying $x\in\pi\left(\left[\nu_{k}\right]\right)$
and \eqref{eq:nu or nuprime is good}. To verify that the corresponding
sum of diameters $\sum_{\nu\in\mathcal{C}}\diam\left(\pi\left(\left[\nu\right]\right)\right)^{b}$
converges, we proceed as follows. If $\nu\notin\mathcal{C}_{\alpha+\epsilon}$
then, by \eqref{eq:nu or nuprime is good}, we can replace $\nu$
by $\nu'\in\mathcal{C}_{\alpha+\epsilon}$, because $\diam\left(\pi\left(\left[\nu\right]\right)\right)\asymp\diam\left(\pi\left(\left[\nu'\right]\right)\right)$
by \eqref{eq:comparable}. This defines a map $\nu\mapsto\nu'$ from
$\mathcal{C}\setminus\mathcal{C}_{\alpha+\epsilon}$ to $\mathcal{C}_{\alpha+\epsilon}$.
Since the involved numbers $\ell_{k}$ and $\ell'_{k}$ in the definition
of $\nu$ and $\nu'$ satisfy \eqref{eq:comparability lk}, we have
that the map $\nu\mapsto\nu'$ is at  most $M$-to-$1$ for some uniform
constant $M\in\N$. Since 
\begin{equation}
\sum_{\omega\in\mathcal{C}_{\alpha+\epsilon}}\diam(\pi([\omega]))^{b}\asymp\sum_{\omega\in\mathcal{C}_{\alpha+\epsilon}}\e^{\left(t(\beta)+\beta(\alpha+\epsilon)+\eta\right)S_{|\omega|}\varphi(\overline{\omega})}\le\sum_{\omega\in\mathcal{C}_{\alpha+\epsilon}}\e^{\left(t(\beta)+\eta\right)S_{\left|\omega\right|}\varphi(\overline{\omega})+\beta S_{\left|\omega\right|}\psi(\overline{\omega})}<\infty,\label{eq:finite s dim hausdorff}
\end{equation}
we therefore conclude that the $b$-dimensional Hausdorff measure
of $\mathcal{F}^{*}(\alpha)$ is finite. Now, first assume that $\alpha\in\left[\alpha_{-},\alpha_{0}\right]$.
Since $\epsilon$ and $\eta$ are arbitrary, it follows that 
\[
\dim_{H}\mathcal{F}^{*}(\alpha)\le\inf_{\beta>0}\left\{ t(\beta)+\beta\alpha\right\} =-t^{*}(-\alpha),
\]
where we have used $\alpha\le\alpha_{0}$ for the last equality. Finally,
if $\alpha<\alpha_{-}$ then $\mathcal{C}_{\alpha+\epsilon}=\emptyset$
if $\alpha+\epsilon<\alpha_{-}$. By the above construction of the
 covering of $\mathcal{F}^{*}(\alpha)$ it thus follows that $\mathcal{F}^{*}(\alpha)=\emptyset$.
The proof is complete.
\end{proof}
\begin{rem}
\label{rem:upper half}For the dimension spectrum of the level sets
$R_{*}(\alpha)$ with $\alpha\ge\alpha_{0}$, the upper bound of the
Hausdorff dimension can also be deduced from \cite[Theorem 1.1 (2)]{MR2864378}
(see also \cite[Proof of Proposition 7.6 ii]{MR3771122}).  
\end{rem}

\newcommand{\etalchar}[1]{$^{#1}$}
\def\cprime{$'$}
\providecommand{\bysame}{\leavevmode\hbox to3em{\hrulefill}\thinspace}
\providecommand{\MR}{\relax\ifhmode\unskip\space\fi MR }
\providecommand{\MRhref}[2]{%
  \href{http://www.ams.org/mathscinet-getitem?mr=#1}{#2}
}
\providecommand{\href}[2]{#2}

\end{document}